\documentclass[10pt,a4paper]{smfart}
\usepackage[latin1]{inputenc}
\usepackage{amsmath}
\usepackage{pstricks}
\usepackage{amsfonts}
\usepackage{amssymb}
\usepackage{amsthm}
\usepackage{lmodern}
\usepackage[all]{xy}
\usepackage{mathrsfs}
\usepackage[english,francais]{babel}
\newtheorem{defi}{Definition}[section]
\newtheorem{thm}[defi]{Theorem}
\newtheorem{prop}[defi]{Proposition}
\newtheorem{lem}[defi]{Lemma}

\newtheorem{coro}[defi]{Corollary}

\newtheorem{rem}[defi]{Remark}
\newtheorem{ex}[defi]{Example}

\author{Jean-Christophe SAN SATURNINO}
\title[Defect of an extension, key polynomials and local uniformization]{Defect of an extension, key polynomials and local uniformization}
\keywords{defect, valued fields, local uniformization, key polynomials}
\subjclass{13A18, 12J10, 12J20, 14B05, 14E15}
\address{Universit\'e Toulouse III Paul Sabatier\\
Institut de Math\'ematiques de Toulouse\\ 
118 route de Narbonne\\
31062 Toulouse cedex 9 (France)}
\email{san@math.ups-tlse.fr}
\urladdr{http://www.math.univ-toulouse.fr/$\sim$san/}

\begin{document}
\selectlanguage{english}
\maketitle
\begin{abstract}
For all simple and finite extension of a valued field, we prove that its defect is the product of the effective degrees of the complete set of key polynomials associated. As a consequence, we obtain a local uniformization theorem for valuations of rank 1 centered on an equicharacteristic quasi-excellent local domain satisfying some inductive assumptions of lack of defect.
\end{abstract}
\setcounter{section}{-1}
\section{Introduction}
\indent In order to obtain a local uniformization theorem, F.-V. Kuhlmann systematized the study of the defect of an extension of valued fields. We know that this defect is a power of $p$, the characteristic of the residual field of the valuation. Another way to approaches the local uniformization problem for a valuation of rank $1$ is studying the set of key polynomials associated to the valuation. The key polynomials have been developped in \cite{vaquie2} and \cite{spivamahboub}, a link of the two approach is given in \cite{mahboub}. When the set of key polynomials do not have limit key polynomials, we proved in \cite{jccar0} that we can obtain a local uniformization theorem for valuation of rank $1$. In particular, we also proved that, if the residual field of the valuation is of characteristic zero, we have local unifomization of the valuation, result known since Zariski.
\\\indent In \cite{vaquie3}, M. Vaqui\'e found a link between the defect and the key polynomials: the jump of the valuation. In the language of \cite{spivamahboub}, this is the degree of a limit key polynomial over the degree of a previous key polynomial whose degree is the same of an infinite family of key polynomials being between this polynomial and the limit key polynomial.
\\\indent In this paper, we want to study more precisely the link between the defect and the key polynomials. We will express the defect as a product of effective degrees and generalize the result of M. Vaqui\'e. The effective degree is closely related with the Newton polygon associated to a key polynomial. We will use this resust to prove a local uniformization theorem for valuations of rank $1$ centered on an equicharacteristic quasi-excellent local domain satisfying some inductive assumptions of defectless of the quotient field.
\\ \\\indent In the first section, we recall some definitions about the center of a valuation and the graded algebras associated to a valuation.
\\\indent In the second section we define the defect of a valuation, roughly speaking, for a valuation $\nu$ of a field $K$ having a unique extension $\mu$ of a field $L$, the defect is the number who gives the equality in this inequality:
$$[L:K]\leq [\Gamma_\mu:\Gamma_\nu][k_\mu:k_\nu],$$
where $\Gamma_\nu$ (resp. $\Gamma_\mu$) is the value group of $\nu$ (resp. $\mu$) and $k_\nu$ (resp. $k_\mu$) is the residual field of $\nu$ (resp. $\mu$). We know that is a power of $p$, the characteristic of $k_\nu$.
\\\indent In the third section we first recall the definitions of key polynomials given in \cite{spivamahboub} or \cite{jccar0}. With this definitions, we give some characterizations of a complete set of key polynomials, this is very useful to prove if a set of polynomials is a complete or a $1$-complete set of key polynomials. For a key polynomial $Q$, we also define the effective degree of a polynomial $P$ corresponding to this key polynomial. When we write the standard expansion of $P$ in terms of $Q$, the greatest power of $Q$ where the valuation of the monomials in $Q$ is maximal is the effective degree of $P$. We prove after that some graded algebras associated to a valuation are euclidean for the effective degree. In the end of this section, we recall that the sequence of effective degrees is decreasing and so, there exist a constant value. If this value is $0$ or $1$, the set of key polynomial considered is $1$-complete.
\\\indent In the fourth section, we first prove that the defect of an extension of a valuation is the product of the effective degree of the limit key polynomials provided by the valuation. The proof uses the notion of key polynomials of W. Mahboub (see \cite{mahboubthese}), so it is true for valuation of any rank. The prinicpal idea is that the Newton plygon has only one face if the field is henselian. We compute after some defects with four examples of M. Vaqui\'e, W. Mahboub and S.D. Cutkosky.
\\\indent In the next section we prove that there is no limit key polynomials for valuation of rank one over a defectless field.
\\\indent In the final section we generalize the results of \cite{jccar0} and prove a local uniformization theorem for a valuation of rank $1$ with some inductive assumptions of local uniformization in lower dimension and no defect for a well defined extension. With this we can prove again the local uniformization for a valuation of rank $1$ in characteristic zero.
\\ \\\indent I want to thank M. Vaqui\'e who suggested me to study the precise link between the defect and the limit key polynomial and B. Teissier for his advice on clarification and structuration of the notion of key polynomials. I also thank W. Mahboub for his ability to produce examples of key polynomial wich verify or not my ideas. Finally, I want to thank M. Spivakovsky for our discussions and his careful reading.
\\ \\\noindent\textbf{Notation.} Let $\nu$ be a valuation of a field $K$. We write $R_\nu=\lbrace f\in K\:\vert\:\nu(f)\geqslant 0\rbrace$, this is a local ring whose maximal ideal is $m_\nu=\lbrace f\in K\:\vert\:\nu(f)> 0\rbrace$. We then denote by $k_\nu=R_\nu/m_\nu$ the residue field of $R_\nu$ and $\Gamma_\nu=\nu(K^*)$.
\\For a field $K$, we will denote by $\overline K$ an algebraic closure of $K$, by $K^{sep}$ a separable closure of $K$ and by $Aut(\overline K\vert K)$ 
 the group of automorphisms of $\overline K\vert K$.
\\If $R$ is a ring and $I$ an ideal of $R$, we will denote by $\widehat{R}^I$ the \textit{ $I$-adic completion} of $R$. When $(R,\mathfrak m)$ is a local ring, we will say the \textit{completion} of $R$ instead of the $\mathfrak m$-adic completion  of $R$ and we will denote it by $\widehat{R}$.
\\For all $P\in Spec(R)$, we note $\kappa(P)=R_P/PR_P$ the residue field of $R_P$.
\\For $\alpha\in\mathbb Z^n$ and $u=(u_1,...,u_n)$ a $n$-uplet of elements of $R$, we write:
\[u^\alpha=u_1^{\alpha_1}...u_n^{\alpha_n}.\]
For $P,Q\in R\left[ X\right] $ with $P=\sum\limits_{i=0}^{n}a_{i}Q^{i}$ and $a_{i}\in R[X]$ such that the degree of $a_{i}$ is strictly less than $Q$, we write:
\[d_{Q}^{\:\circ}(P)=n.\]
If $Q=X$, we will note simply $d^{\:\circ}(P)$ instead of $d_{X}^{\:\circ}(P)$.
\\Finally, if $R$ is a domain, we denote by $Frac(R)$ its quotient field.

\section{Center of a valuation, graded algebra associated and saturation}
In this section, $\Gamma$ will be a totally ordered commutative group.
\begin{defi}\label{centreval}
Let $R$ be a ring and $P$ a prime ideal. A valuation $\nu:R\rightarrow\Gamma\cup\{\infty\}$ \textbf{centered on} $\boldsymbol P$ is the data of a minimal prime ideal $P_\infty$ of $R$ included in $P$ and a valuation of the quotient field of $R/P_\infty$ centered on $P/P_\infty$. The ideal $P_\infty$ is the support of the valuation, that is $P_\infty=\nu^{-1}(\infty)$.
\\If $R$ is a local ring with maximal ideal $\mathfrak{m}$, we will say that $\nu$ is \textbf{centered on} $\boldsymbol R$ to say that $\nu$ is centered on $\mathfrak{m}$.
\end{defi}

\begin{defi}\label{alggrad}
Let $R$ be a ring and $\nu:R\rightarrow \Gamma\cup\lbrace\infty\rbrace$ a valuation centered on a prime ideal of $R$. For all $\alpha\in\nu(R\setminus\lbrace 0\rbrace)$, we define the ideals:
\[\mathcal P_\alpha=\lbrace f\in R\:\vert\:\nu(f)\geqslant\alpha\rbrace;\]
\[\mathcal P_{\alpha,+}=\lbrace f\in R\:\vert\:\nu(f)>\alpha\rbrace.\]
Then we define \textbf{the graded algebra of} $\boldsymbol{R}$ \textbf{associated to} $\boldsymbol{\nu}$ by:
\[gr_{\nu}(R)=\bigoplus\limits_{\alpha\in\nu(R\setminus\lbrace 0\rbrace)}P_{\alpha}/P_{\alpha,+}.\]
The algebra $gr_\nu(R)$ is a domain. For $f\in R\setminus\lbrace 0\rbrace$, we write $in_\nu(f)$ its image in $gr_\nu(R)$.
\end{defi}
If $R$ is a local domain, we define another graded algebra by:
\begin{defi}
Let $R$ be a local domain, $K=Frac(R)$ and $\nu:K^*\twoheadrightarrow \Gamma\cup\lbrace\infty\rbrace$ a valuation of $K$ centered on $R$. For all $\alpha\in\Gamma$, we define the $R_\nu$-submodules of $K$ following:
\[P_\alpha=\lbrace f\in K\:\vert\:\nu(f)\geqslant\alpha\rbrace\cup\{0\};\]
\[P_{\alpha,+}=\lbrace f\in K\:\vert\:\nu(f)>\alpha\rbrace\cup\{0\}.\]
We define then \textbf{the graded algebra associated at} $\boldsymbol{\nu}$ by:
\[G_{\nu}=\bigoplus\limits_{\alpha\in\Gamma}P_{\alpha}/P_{\alpha,+}.\]
For $f\in K^*$, we write $in_\nu(f)$ his image in $G_\nu$.
\end{defi}
\begin{rem}
\textup{We have the natural injective function:
\[gr_\nu(R)\hookrightarrow G_\nu.\]}
\end{rem}
\begin{defi}\label{defisature}
Let $G$ be a graded algebra without zero divisors. We called \textbf{saturation of} $\boldsymbol G$ the graded agebra $G^*$ define by:
\[G^{*}=\left\lbrace \left.\dfrac{f}{g}\:\right|\:f,g\in G,\:g\textit{ homog\`ene},\:g\neq 0\right\rbrace.\]
We say that $G$ is \textbf{saturated} if $G=G^*$.
\end{defi}
\begin{rem}
\textup{For all graded algebra $G$, we have:
\[G^*=\left(G^*\right)^*.\]
In particular, $G^*$ is always saturated.}
\end{rem}
\begin{ex}
\textup{Let $\nu$ be a valuation centered on a local ring $R$. Then:
\[G_\nu=\left(gr_\nu(R)\right)^*.\]
In particular, $G_\nu$ is saturated.}
\end{ex}

\section{Defectless fields}
\begin{defi}(\cite{vaquie3}, D\'efinition 1.2)
A valued field $(K,\nu)$ is said to be \textbf{henselian} if for all algebraic extension $L\vert K$, there exsits a unique valuation $\mu$ of $L$ which extends $\nu$.
\\We call \textbf{henselization of} $\boldsymbol{(K,\nu)}$ all extension $(K^{h},\mu)$ of $(K,\nu)$ such that $(K^h,\mu)$ is henselian and for all henselian valued field $(L,\nu')$ and all immersion $\sigma:(K,\nu)\hookrightarrow (L,\nu')$, there exists a unique immersion $\sigma': (K^h,\mu)\hookrightarrow (L,\nu')$ which extends $\sigma$.
\[\xymatrix{
    (K^h,\mu) \ar@{.>}[r]^-{\sigma'}   & (L,\nu')  \\
    (K,\nu)\ar[ru]_-{\sigma} \ar[u]&
  }\]
 \rput(6.14,1.45){$\circlearrowleft$}
\end{defi}
\begin{rem}
\textup{All the henselizations of a given valued field are isomorphic, it is its smaller henselian extensions. Moreover, the henselization $(K^h,\mu)$ of $(K,\nu)$ is an immediate extension, that is $\Gamma_\nu=\Gamma_\mu$ and $k_\nu=k_\mu$ (for a proof, we can see Theorem 7.42 of \cite{kulhlivre}).}
\end{rem}
The henselization of a given valued field $(K,\nu)$ can be constructed explicitly, we give here the construction proposed by Kuhlmann in \cite{kulhlivre}, Chapter 7. Consider $\mu$ to be an extension of $\nu$ on $\overline K$. Write:
\[G^d(\overline K\vert K,\mu)=\lbrace \sigma\in Aut(\overline K\vert K)\:\vert\:\forall\:\alpha\in\overline K,\:\mu(\sigma(\alpha))=\mu(\alpha)\rbrace,\]
\[K^{h(\mu)}=\lbrace \alpha\in K^{sep}\:\vert\:\forall\:\sigma\in G^d(\overline K\vert K,\mu), \sigma(\alpha)=\alpha\rbrace.\]
$(K^{h(\mu)},\mu)$ is an henselization of $(K,\nu)$, it is an algebraic separable immediate extension of $(K,\nu)$.
\begin{defi}\label{defidefautlocal}
Let $K$ be a field with a valuation $\nu$. Consider $L$ a finite extension of $K$ and write $\mu_1,...,\mu_g$ the valuations of $L$ which extend $\nu$. For $1\leqslant i\leqslant g$, choose a valuation $\overline\mu_i$ of $\overline K$ whose restriction on $L$ is $\mu_i$. We write $K^{h(\overline\mu_i)}$ the field construct previously, it is an henselization of $(K,\nu)$. Moreover, it is a subfield of $L^{h(\overline\mu_i)}:=L.K^{h(\overline\mu_i)}$ which is himself an henselization of $(L,\mu_i)$. We define the \textbf{defect of the extension} $\boldsymbol{L\vert K}$ \textbf{in} $\boldsymbol{\mu_i}$ by:
\[d_{L\vert K}(\mu_i,\nu)=\dfrac{\left[L^{h(\overline\mu_i)}:K^{h(\overline\mu_i)}\right]}{e_if_i},\]
where $e_i=\left[\Gamma_{\mu_i}:\Gamma_\nu\right]$ and $f_i=\left[k_{\mu_i}:k_\nu\right]$.
\end{defi}
\begin{rem}\label{ostrow}\textup{By a theorem of Ostrowski (Lemma 11.17 of \cite{kulhlivre}), we can show that $d_{L\vert K}(\mu_i,\nu)=p^{a_i}$ where $a_i\in\mathbb N$ and $p=char(k_\nu)$.}
\end{rem}
\begin{prop}(\cite{kulhlivre}, Lemma 7.46)\label{thmdefautdegre}
Let $K$ be a field with a valuation $\nu$. Consider $L$ a finite extension of $K$ and write $\mu_1,...,\mu_g$ the valuations of $L$ which extend $\nu$. Then:
\[\left[L:K\right]=\sum\limits_{i=1}^gd_ie_if_i,\]
where $d_i=d_{L\vert K}(\mu_i,\nu)$, $e_i=\left[\Gamma_{\mu_i}:\Gamma_\nu\right]$ and $f_i=\left[k_{\mu_i}:k_\nu\right]$.
\end{prop}
\begin{defi}
Let $K$ be a field with a valuation $\nu$. Consider $L$ a finite extension of $K$ and write $\mu_1,...,\mu_g$ the valuation of $L$ which extend $\nu$. We call \textbf{global defect of the extension} $\boldsymbol{L\vert K}$ \textbf{on} $\boldsymbol{\nu}$ the quotient:
\[d_{L\vert K}(\nu)=\dfrac{\left[L:K\right]}{\sum\limits_{i=1}^ge_if_i},\]
where $e_i=\left[\Gamma_{\mu_i}:\Gamma_\nu\right]$ et $f_i=\left[k_{\mu_i}:k_\nu\right]$.
\\We say that $L$ is a \textbf{defectless extension} of $K$ if:
\[d_{L\vert K}(\nu)=1.\]
\end{defi}
\begin{rem}
\textup{If the extension $L\vert K$ is normal, then the global defect is a power of $p$, where $p=char(k_{\nu})$, because $d_{L\vert K}(\nu)=d_{L\vert K}(\mu_i,\nu)$, for all $i\in\{1,...,g\}$ (see \cite{kulhlivre}, Lemma 11.3) ; otherwise it is a rational number.}
\end{rem}
\begin{rem}
\textup{$d_{L\vert K}(\nu)=1 \Leftrightarrow \forall\:i\in\lbrace 1,...,g\rbrace,\:d_{L\vert K}(\mu_i,\nu)=1.$ }
\end{rem}
\begin{defi}
A field $K$ is said to be \textbf{defectless} if all finite extension of $K$ is defectless.
\end{defi}
\begin{prop}\label{cropasdefaut}(\cite{kulhlivre}, Theorem 11.23)
All field with a valuation $\nu$ such that $char(k_\nu)=0$ is a defectless field.
\end{prop}
\section{Key polynomials and effective degree}
Let $K\hookrightarrow K(x)$ be a simple transcendental extension of fields. Let $\mu'$ be a valuation of $K(x)$, write $\mu:=\mu'_{\vert\:K}$. Write $G$ the group of values of $\mu'$ and $G_{1}$ the smallest non-zero isolated subgroup of $G$. Suppose that the rank of $\mu$ is $1$ and $\mu'(x)>0$. Finally, for $\beta\in G$, recall that:
\[P'_{\beta}=\lbrace f\in K(x)\:\vert\:\mu'(f)\geqslant\beta\rbrace\cup\lbrace 0\rbrace;\]
\[P'_{\beta,+}=\lbrace f\in K(x)\:\vert\:\mu'(f)>\beta\rbrace\cup\lbrace 0\rbrace;\]
\[G_{\mu'}=\bigoplus\limits_{\beta\in G}P'_{\beta}/P'_{\beta,+};\]
and $in_{\mu'}(f)$ is the image of $f\in K(x)$ in $G_{\mu'}$.
\begin{defi}
A \textbf{complete set of key polynomials} for $\mu'$ is a well ordered collection:
\[\textbf{Q}=\lbrace Q_{i}\rbrace_{i\in\Lambda}\subset K[x]\]
such that, for all $\beta\in G$, the additive group $P'_{\beta}\cap K[x]$ is generated by the products of the form $a\prod\limits_{j=1}^{s}Q_{i_{j}}^{\gamma_{j}}$, $a\in K$, such that $\sum\limits_{j=1}^{s}\gamma_{j}\mu'\left(Q_{i_{j}}\right)+\mu(a)\geqslant\beta$.
\\The set is said to be \textbf{1-complete} if the precedent condition occurs for all $\beta\in G_1$ and if for all $i\in\Lambda$, $\mu'(Q_i)\in G_1$.
\end{defi}
\begin{thm}(\cite{spivamahboub}, Theorem 7.11)\label{existpolycle}
There exists a collection $\textbf{Q}=\lbrace Q_{i}\rbrace_{i\in\Lambda}$ which is a $1$-complete set of key polynomials.
\end{thm}
By the Theorem \ref{existpolycle}, we know that there exists a $1$-complete set of key polynomials $\textbf{Q}=\lbrace Q_{i}\rbrace_{i\in\Lambda}$ and that the order type of $\Lambda$ is at most $\omega\times\omega$. If $K$ is defectless, we will see that the order type of $\Lambda$ is at most $\omega$ at so, there is no limit key polynomials, this is particulary the case if $char(k_\mu)=0$. For all $i\in\Lambda$, write $\beta_i=\mu'(Q_i)$.

Let $l\in\Lambda$, write:
\[\alpha_{i}=d_{Q_{i-1}}^{\:\circ}(Q_{i}),\:\forall\:i\leqslant l;\]
\[\boldsymbol{\alpha_{l+1}}=\lbrace \alpha_{i}\rbrace_{i\leqslant \:l};\]
\[\textbf{Q}_{l+1}=\lbrace Q_{i}\rbrace_{i\leqslant \:l}.\]
We also use the notation $\overline{\gamma}_{l+1}=\lbrace \gamma_{i}\rbrace_{i\leqslant\:l}$ where the $\gamma_{i}$ are all zero except a finite number, $\textbf{Q}_{l+1}^{\overline{\gamma}_{l+1}}=\prod\limits_{i\leqslant\: l}Q_{i}^{\gamma_{i}}$.

\begin{defi}
A multiindex $\overline{\gamma}_{l+1}$ is said to be \textbf{standard with respect to} $\boldsymbol{\alpha_{l+1}}$ if $0\leqslant \gamma_{i}<\alpha_{i+1}$, for $i\leqslant l$.
\\ An \textbf{l-standard monomial in} $\boldsymbol{Q_{l+1}}$ is a product of the form $c_{\overline{\gamma}_{l+1}}\textbf{Q}_{l+1}^{\overline{\gamma}_{l+1}}$, where $c_{\overline{\gamma}_{l+1}}\in K$ and $\overline{\gamma}_{l+1}$ is standard with respect to $\boldsymbol{\alpha_{l+1}}$.
\\ An \textbf{l-standard expansion not involving} $\boldsymbol{Q_{l}}$ is a finite sum $\sum\limits_{\beta}S_{\beta}$ of $l$-standard monomials  not involving $Q_{l}$, where $\beta$ ranges over a certain finite subset of $G_{+}$ and $S_{\beta}=\sum\limits_{j} d_{\beta,j}$ is a sum of standard monomials of value $\beta$ such that $\sum\limits_{j} in_{\mu'}(d_{\beta,j})\neq 0$.
\end{defi}
\begin{defi}
Let $f\in K[x]$ and $i\leqslant l$, an \textbf{i-standard expansion  of f} is an expression of the form:
\[f=\sum\limits_{j=0}^{s_{i}}c_{j,i}Q_{i}^{j},\]
where $c_{j,i}$ is an $i$-standard expansion not involving $Q_{i}$.
\end{defi}
\begin{rem}\textup{
Such expansion exists by Euclidean division and is unique in the sense of that the $c_{j,i}\in K[x]$ are unique. More precisely, if $i\in\mathbb N$, we can prove by induction that the $i$-standard expansion is unique.}
\end{rem}
\begin{defi}
Let $f\in K[x]$, $i\leqslant l$ and $f=\sum\limits_{j=0}^{s_{i}}c_{j,i}Q_{i}^{j}$ an $i$-standard expansion of $f$. We define the \textbf{i-troncation of} $\boldsymbol{\mu'}$, denoted by $\mu_{i}'$, as being the pseudo-valuation:
\[\mu_{i}'(f)=\min_{0\leqslant j\leqslant s_{i}}\lbrace j\mu'(Q_{i})+\mu'(c_{j,i})\rbrace.\]
\end{defi}
\begin{rem}\label{inegalitetronque}\textup{
We can prove that $\mu_{i}'$ is a valuation. Moreover, we have:
\[\forall \:f\in K[x],\: i\in\Lambda, \:\mu_{i}'(f)\leqslant \mu'(f).\]}
\end{rem}

\begin{prop}
Let $\lbrace Q_{i}\rbrace_{i\in\Lambda}\subset K\left[ x\right]$. Write $H=\lbrace f\in K[x] \:\vert\:\mu'(f)\notin G_1\rbrace$ and, for a graded algebra $G$, denote by $G^*$ his saturation (see Definition \ref{defisature}). Consider the following assertions:
\begin{enumerate}
\item The set $\lbrace Q_{i}\rbrace_{i\in\Lambda}$ is $1$-complete;
\item $\forall\:f\in K[x]\setminus H$, $\exists\:i\in\Lambda,\:\mu_i'(f)=\mu'(f)$.
\item For all $i\in\Lambda$, $Q_{i}\notin H$ and there exists $h\in H$ monic or zero such that the set $\lbrace Q_{i}\rbrace_{i\in\Lambda}\cup\{h\}$ is complete;
\item There exists $h\in H$ monic or zero such that $G_{\mu'}=\left( G_\mu\left[\{in_{\mu'}(Q_i)\}_{i\in\Lambda},in_{\mu'}(h)\right]\right)^*$.
\end{enumerate}
Then \textup{1.} $\Leftrightarrow$ \textup{2.} $\Leftrightarrow$ \textup{3.} and \textup{1.} $\Rightarrow$ \textup{4.}\:.
\end{prop}
\noindent\textit{Proof}: Show that $\lbrace Q_{i}\rbrace_{i\in\Lambda}$ is $1$-completed if and only if for all $f\in K[x]\setminus H$ there exists $i\in\Lambda$, $\mu_i'(f)=\mu'(f)$. \\Suppose that $\lbrace Q_{i}\rbrace_{i\in\Lambda}$ is $1$-completed and consider $f\in K[x]$ such that $\mu'(f)\in G_1$. Write $\beta=\mu'(f)$, then $P'_{\beta}\cap K[x]$ is generated by the products of the form $a\prod\limits_{j=1}^{s}Q_{i_{j}}^{\gamma_{j}}$, $a\in K$, with $\sum\limits_{j=1}^{s}\gamma_{j}\mu'\left(Q_{i_{j}}\right)+\mu(a)\geqslant\beta=\mu'(f)$, thus $\mu'_{i_s}(f)\geqslant \mu'(f)$. Moreover, by definition, for all $i\in\Lambda$, $\mu_i'(f)\leqslant \mu'(f)$. We conclude that $\mu_{i_s}'(f)=\mu'(f)$. Reciprocally, if all $f\in K[x]$ such that $\mu'(f)\in G_1$ is as $\mu_i'(f)=\mu'(f)$, for a  $i\in\Lambda$, then, $f$ writes such a sum of monomials in $\boldsymbol{Q}_{i+1}$ of value at least $\mu'(f)\geqslant \beta$. So $f$ is in the ideal generated by all this monomials.
We can show in the same way that $\lbrace Q_{i}\rbrace_{i\in\Lambda}$ is completed if and only if for all $f\in K[x]$ and for all $i\in\Lambda$, $\mu_i'(f)=\mu'(f)$.
\\\indent Show that $\lbrace Q_{i}\rbrace_{i\in\Lambda}$ is $1$-completed if and only if for all $i\in\Lambda$, $Q_{i}\notin H$ and there exists $h\in H$ monic or zero such that the set $\lbrace Q_{i}\rbrace_{i\in\Lambda}\cup\{h\}$ is completed. \\If $\lbrace Q_{i}\rbrace_{i\in\Lambda}$ is completed then, by definition, for all $i\in\Lambda$, $\mu'(Q_i)\in G_1$ and take $h=0$. If $\lbrace Q_{i}\rbrace_{i\in\Lambda}$ is not completed, there exists an element $h\in K[x]$ such that $\mu'(h)\notin G_1$ and $\mu_i'(h)<\mu'(h)$, for all $i\in \Lambda$. Take $h$ of minimal degree and monic satifying the previous inequality. For $f\in K[x]$, consider the standard expansion in terms of $h$: $f=\sum\limits_{j=0}^{s}c_{j}h^j$. Write $m=\min\{j\in\lbrace 1,...,s\rbrace\:\vert\: c_j\neq 0\}$. Denote by $\mu_h'$ the troncated valuation in terms of $h$, since $\mu'(h)\notin G_1$ then:
\[\mu_h'(f)=\mu_h'(c_mh^m)=\mu'(c_mh^m)=\mu'(f).\]
The reciproque is obvious.
\\\indent Finally, if $\lbrace Q_{i}\rbrace_{i\in\Lambda}$ is $1$-completed then there exists $h\in H$ monic or zero such that the set $\lbrace Q_{i}\rbrace_{i\in\Lambda}\cup\{h\}$ is completed, that is $P'_\beta\cap K[x]$ is generated by elements of the form $a\prod\limits_{j=1}^sg_j^{\gamma_j}$ where $g_j\in\lbrace Q_{i}\rbrace_{i\in\Lambda}\cup\{h\}$. So we deduce that $G_{\mu'}=\left( G_\mu\left[\{in_{\mu'}(Q_i)\}_{i\in\Lambda},in_{\mu'}(h)\right]\right)^*$.\\\qed

\begin{rem}
\textup{Note that 4. $\not\Rightarrow$ 1., consider the example 2.2 of \cite{ideimp}:
\\Let $k$ a field and $\iota:k(u,v)\hookrightarrow k\left[\left[t\right]\right]$ a monomorphism defined by $\iota(v)=t$ and $\iota(u)=\sum\limits_{i\geqslant 1}c_it^i$ where $c_i\in k^*$ and $Q_\infty=u-\sum\limits_{i\geqslant 1}c_iv^i$ is transcendantal over $k(u,v)$. Denote by $\mu'$ the valuation of $k(u,v)$ induced by the $t$-adic valuation of $k\left[\left[t\right]\right]$ through $\iota$ and $\mu$ by the restriction of $\mu'$ on $K=k(v)$. The set of key polynomials associated to the extension $K\hookrightarrow K(u)$ is $\{Q_i\}_{i\geqslant 1}$ with $Q_1=u$ and $Q_i=u-\sum\limits_{j= 1}^{i-1}c_iv^i$ while $G_{\mu'}=k_\mu[in_{\mu'}(u)]$.}
\end{rem}

\indent For the next definition, we use the terminology of \cite{vaquie5} following the notations of \cite{spivamahboub}.
\begin{defi}
Let $h\in K[x]$, consider its $i$-standard expansion $h=\sum\limits_{j=0}^{s_{i}}c_{j,i}Q_{i}^j$. We call by the $\boldsymbol{i}$\textbf{-th effective degree of h} the natural number:
\[\delta_i(h)=\max \lbrace j\in\lbrace 0,...,s_{i}\rbrace\:\vert\:j\beta_i+\mu'\left(c_{j,i}\right)=\mu'_i(h)\rbrace.\]
By convention, $\delta_i(0)=-\infty$.
\end{defi}
\begin{rem}\label{remdetadegre}
\textup{If we denote by: $$in_i(h)=\sum\limits_{j\in S_i(h,\beta_i)}in_{\mu'}(c_{j,i})X^j\in G_{\mu}[in_{\mu'}(\textbf Q_i),X]$$
where $S_i(h,\beta_i)=\lbrace j\in\lbrace 0,...,s_{i}\rbrace\:\vert\:j\beta_i+\mu'\left(c_{j,i}\right)=\mu'_i(h)\rbrace$, then $\delta_i(h)=d^{\:\circ}_X(in_i(h))$.
\\Moreover, if $in_{\mu_i'}(h)=in_{\mu_i'}(g)$ then $\delta_i(h)=\delta_i(g)$. The definition of the effective degree  extends naturally to $gr_{\mu_i'}(K[x])$.}
\end{rem}
\begin{rem}
\textup{Remind that, by Proposition 5.2 of \cite{spivamahboub}, for $l\in\Lambda$ an ordinal number, the sequence $(\delta_{l+i}(h))_{i\in\mathbb N^*}$ decreases. Thus there exists $i_0\in\mathbb N^*$ such that $\delta_{l+i_0}(h)=\delta_{l+i_0+i}(h)$, for all $i\geqslant 1$ and we denote this common value by $\delta_{l+\omega}(h)$ or $\delta_{l+\omega}$ if no confusion is possible. Note that $\delta_{l+\omega}$ may to be $0$.}
\end{rem}
The next three lemmas generalize the lemmas 2.12, 2.13 and 2.14 of \cite{favrejon}, the proofs are quasi the same.
\begin{lem}\label{deltaprod}
Let $l\in\Lambda$ and $f,g\in K[x]$:
\[\delta_{l}(fg)=\delta_{l}(f)+\delta_{l}(g).\] 
\end{lem}
\noindent\textit{Proof}: It is an immediate consequence of the Remark \ref{remdetadegre}.\\\qed
\begin{lem}\label{deltainvzero}
For all $h\in K[x]$, $\delta_{i}(h)=0$ if and only if $in_{\mu_i'}(h)$ is a unit of $gr_{\mu_i'}(K[x])$.
\end{lem}
\noindent\textit{Proof}: If $\delta_{i}(h)=0$ then $\mu_i'(h)=\mu'(h)$ and $in_{\mu_i'}(h)=in_{\mu_i'}(c_{0,i})$ in $gr_{\mu_i'}(K[x])$. As the polynomial $Q_i$ is irreducible in $K[x]$ and $d^{\:\circ}(c_{0,i})<d^{\:\circ}(Q_i)$, the polynomials $Q_i$ and $c_{0,i}$ are coprime in $K[x]$. Then, there exists $U,V\in K[x]$ such that $UQ_i+Vc_{0,i}=1$. We deduce that $\mu_i'(Vc_{0,i})=\mu_i'(1)<\mu_i'(UQ_i)$ and $h$ is a unit of $gr_{\mu_i'}(K[x])$.
\\Reciprocally, it is sufficient to apply the Lemma \ref{deltaprod} and the Remark \ref{remdetadegre}.\\\qed
\begin{lem}
For all $i\in\Lambda$, the ring $gr_{\mu_i'}(K[x])$ is euclidean for $\delta_i$, ie:
$$\forall g,h\in K[x],\:h\neq 0,\:\exists\:Q,R\in K[x],\:g=hQ+R\textit{ in }gr_{\mu_i'}(K[x])\textit{ and }0\leqslant \delta_{i}(R)<\delta_i(h).$$
\end{lem}
\noindent\textit{Proof}: Write $h=\sum\limits_{j=0}^{s_{i}}c_{j,i}Q_{i}^j$, we can suppose without loss of generalities then $c_{j,i}=0$ for $j>\delta_i(h)$ because $\mu'(c_{j,i}Q_{i}^j)>\mu_i'(c_{j,i}Q_{i}^j)$ and so $in_{\mu_i'}\left(\sum\limits_{j=\delta_i(h)+1}^{s_{i}}c_{j,i}Q_{i}^j\right)=0$. We can also suppose that $c_{\delta_i(h),i}=1$ by Lemma \ref{deltainvzero}. Note that $d^{\:\circ}(h)=\delta_i(h). d^{\:\circ}(Q_i)$. By euclidean division in $K[x]$, there exists $P,Q\in K[x]$ such that $g=hQ+P$ with $0\leqslant d^{\:\circ}(P)<d^{\:\circ}(h)$. Let $P=\sum\limits_{j=0}^{t_{i}}p_{j,i}Q_{i}^j$ the $i$-standard expansion of $P$. If we write $R=\sum\limits_{j=0}^{\delta_{i}(P)}p_{j,i}Q_{i}^j$, we obtain that $h=gQ+R$ in $gr_{\mu_i'}(K[x])$ and:
\[\delta_{i}(P). d^{\:\circ}(Q_i)\leqslant d^{\:\circ}(p_{\delta_{i}(P),i})+\delta_{i}(P). d^{\:\circ}(Q_i)=d^{\:\circ}(R)\leqslant d^{\:\circ}(P)<d^{\:\circ}(h)=\delta_i(h).d^{\:\circ}(Q_i).\]
Thus, $\delta_{i}(R)=\delta_{i}(P)<\delta_i(h)$.\\\qed

\indent The construction of the key polynomials is recursive (see \cite{mahboub}, \cite{spivamahboub} and \cite{spivaherrera}). For $l\in\mathbb N^*$, we construct a set of key polynomials $\textbf{Q}_{l+1}=\lbrace Q_{i}\rbrace_{1\leqslant i\leqslant \:l}$; there are two cases:
\begin{enumerate}
\item[(1)] $\exists\:l_0\in\mathbb N,\:\beta_{l_0}\notin G_1$;
\item[(2)] $\forall\:l\in\mathbb N,\:\beta_l\in G_1$.
\end{enumerate}
In case (1), we stop the construction; the set $\textbf{Q}_{l_0}=\lbrace Q_{i}\rbrace_{1\leqslant i\leqslant \:l_0-1}$ is by definition a $1$-complete set of key polynomials and $\Lambda=\lbrace 1,...,l_0-1\rbrace$. Note that the set $\textbf{Q}_{l_0+1}$ is a complete set of key polynomials.
\\In case (2), if $K$ is defectless, we will prove in Proposition \ref{polycleencar0} that the set $\textbf{Q}_{\omega}=\lbrace Q_{i}\rbrace_{i\geqslant 1}$ is infinite and $\Lambda=\mathbb N^*$ ; the next propositions will ensure us that the set of key polynomials obtained is also $1$-complete.
\\\indent The next lemma, very useful, allows us to note that there no exists any increasing bounded sequence with valuations of rank $1$.
\begin{lem}\label{pasdesuitecroissantebornee}
Let $\nu$ be a valuation of rank $1$ centered on a local noetherian ring $R$. Denote by $P_\infty$ the support of $\nu$. Then, $\nu\left( R\setminus P_\infty\right)$ does not contain any increasing bounded infinite sequence.
\end{lem}
\begin{rem}
\textup{Recall the Definition \ref{centreval}, the data of a valuation $\nu$ centered on a local ring $(R,\mathfrak{m})$ is the data of a minimal prime ideal  $P_\infty$ of $R$ (the support of the valuation) and of a valuation $\nu'$ of the quotient filed of $R/P_\infty$ such that $R/P_\infty\subset R_{\nu'}$ and $\mathfrak m/P_\infty=R/P_\infty\cap m_{\nu'}$.}
\end{rem}
\noindent\textit{Proof}: Let $\left(\beta_i\right)_{i\geqslant 1}$ an increasing infinite sequence of $\nu\left( R\setminus P_\infty\right)$ bounded by $\beta$. This sequence corresponds to a decreasing infinite sequence of ideals of $R/P_\beta$. It is sufficient to prove that $R/P_\beta$ have finite length. Write $\mathfrak m$ the maximal ideal of $R$, $\nu(\mathfrak m)=\min\left\lbrace\nu\left( R\setminus P_\infty\right)\setminus\lbrace 0\rbrace\right\rbrace$ and $\Gamma$ the value group of $\nu$. Note that the group $\nu\left( R\setminus P_\infty\right)$ is archimedean. Indeed, by absurde, if $\nu\left( R\setminus P_\infty\right)$ is not archimedean, there exists $\alpha,\:\beta\in\nu\left( R\setminus P_\infty\right)$, $\beta\neq 0$ such that, for all $n\geqslant 1$, $n\beta\leqslant\alpha$. Particularly, the set:
\[\lbrace \gamma\in\Gamma\:\vert\:\exists\:n\in\mathbb{N}\setminus\lbrace 0\rbrace,\:-n\beta<\gamma<n\beta\rbrace\]
is a non trivial isolated subgroup of $\Gamma$.
\\We deduce that ther exists $n\in\mathbb N$ such that:
\[\beta\leqslant n\nu(\mathfrak m).\]
Thus, $\mathfrak m^n\subset P_\beta$ and therefore, there exists a surjective map:
\[R/\mathfrak m^n\twoheadrightarrow R/P_\beta.\]
\qed

\begin{prop}\label{sibornealorscomplet}(\cite{spivamahboub}, Proposition 3.30) Suppose that we have construct an infinite set of key polynomials $\textbf{Q}_{\omega}=\lbrace Q_{i}\rbrace_{i\geqslant 1}$ such that, fo all $i\in\mathbb N^*$, $\beta_i\in G_1$. Suppose further that the sequence $\lbrace\beta_i\rbrace_{i\geqslant 1}$ is not bounded in $G_1$. Then, the set of key polynomials $\textbf{Q}_{\omega}$ is $1$-complete.
\end{prop}
\noindent\textit{Proof}: It is sufficient to prove that, for all $\beta\in G_1$ and for all $h\in K\left[x\right]$ such that $\mu'(h)=\beta$, $h$ belongs to the $R_\mu$-submodule of $K\left[x\right]$ generated by all the monomials of the form $a\prod\limits_{j=1}^{s}Q_{i_{j}}^{\gamma_{j}}$, $a\in K$, such that $\mu'\left(a\prod\limits_{j=1}^{s}Q_{i_{j}}^{\gamma_{j}}\right)\geqslant\beta$.
\\Therefore consider $h\in K\left[x\right]$ such that $\mu'(h)\in G_1$. Write $h=\sum\limits_{j=0}^dh_jx^j$, we can suppose, without loss of generalities, that:
\[\forall\:j\in\lbrace 0,...,d\rbrace,\:\mu(h_j)\geqslant 0.\]
Indeed, otherwise it is sufficient to multiply $h$ by an element of $K$ appropriately chosen.
\\Since the sequence $\lbrace\beta_i\rbrace_{i\geqslant 1}$ is not bounded in $G_1$, there exists $i_0\in\mathbb N^*$ such that:
\[\mu'(h)<\beta_{i_0}.\]
Then denote by $h=\sum\limits_{j=0}^{s_{i_0}}c_{j,i_0}Q_{i_0}^j$, the  $i_0$-standard expansion of $h$. This expansion is obtain by euclidean division, seen the choice made on the coefficients of $h$ and, since the sequence $\left\lbrace\frac{\beta_i}{d^{\:\circ}\left(Q_i\right)}\right\rbrace_{i\geqslant 1}$ is strictly increasing (it is sufficient to take the $(i-1)$-standard expansion of $Q_i$), we prove easily that:
\[\forall\:j\in\lbrace 0,...,s_{i_0}\rbrace,\:\mu\left(c_{j,i_0}\right)\geqslant 0.\]
Recall that, by construction of the key polynomials, for $j\in\lbrace 0,...,s_{i_0}\rbrace$, $\mu'_{i_0}\left(c_{j,i_0}\right)=\mu'\left(c_{j,i_0}\right)$. We then deduce that:
\[\forall\:j\in\lbrace 1,...,s_{i_0}\rbrace,\:\mu'\left(c_{j,i_0}Q_{i_0}^j\right)=\mu'_{i_0}\left(c_{j,i_0}Q_{i_0}^j\right)>\mu'(h).\]
Thus, $\mu'(h)=\mu'\left(c_{0,i_0}\right)$ and so, $h$ is a sum of monomials in $\textbf{Q}_{i_0+1}$ of valuation at least $\mu'(h)$ (in particular, $\mu'_{i_0}(h)=\mu'(h)$).\\\qed
\\ \\\indent 
We consider now two cases:
\begin{enumerate}
\item[(1)] $\sharp\lbrace i\geqslant 1\:\vert\:\alpha_i>1\rbrace=+\infty$;
\item[(2)] $\sharp\lbrace i\geqslant 1\:\vert\:\alpha_i>1\rbrace<+\infty$.
\end{enumerate}
In case (1), with the Proposition \ref{sialphainfini}, we can prove that the infinite set of key polynomials is always $1$-complete, independently of the characteristic of $k_\mu$. In case (2), if the effective degree is always $1$ and if the set of key polynomials $\textbf{Q}_{\omega}=\lbrace Q_{i}\rbrace_{i\geqslant 1}$ is not complete, we will prove in the Proposition \ref{sialphafini} that the sequence $\lbrace\beta_i\rbrace_{i\geqslant 1}$ is never bounded. In this case, with the Proposition \ref{sibornealorscomplet}, we deduce that the set of key polynomials $\textbf{Q}_{\omega}=\lbrace Q_{i}\rbrace_{i\geqslant 1}$ is also $1$-complete.
\begin{prop}\label{sialphainfini}(\cite{spivamahboub}, Corollary 5.8)
Suppose that we have construct an infinite set of key polynomials $\textbf{Q}_{\omega}=\lbrace Q_{i}\rbrace_{i\geqslant 1}$ such that, for all $i\in\mathbb N^*$, $\beta_i\in G_1$. More, suppose that the set $\lbrace i\geqslant 1\:\vert\:\alpha_i>1\rbrace$ is infinite. Then, $\textbf{Q}_{\omega}$ is a $1$-complete set of key polynomials.
\end{prop}
\noindent\textit{Proof}: Let $h\in K[x]$, as in the proof of the Proposition \ref{sibornealorscomplet}, it is sufficient to show that $\mu'_{i}(h)=\mu'(h)$ for a $i\geqslant 1$. But, if we write:
\[\delta_i(h)=\max S_i(h,\beta_i),\]
where:
\[S_i(h,\beta_i)=\lbrace j\in\lbrace 0,...,s_{i}\rbrace\:\vert\:j\beta_i+\mu'\left(c_{j,i}\right)=\mu'_i(h)\rbrace,\]
\[h=\sum\limits_{j=0}^{s_{i}}c_{j,i}Q_{i}^j,\]
by (1) of Proposition 5.2 of \cite{spivamahboub}, we have:
\[\alpha_{i+1}\delta_{i+1}(h)\leqslant\delta_i(h),\:\forall\:i\geqslant 1.\]
We deduce that if $\delta_i(h)>0$ and $\alpha_{i+1}>1$:
\[\delta_{i+1}(h)<\delta_i(h),\forall\:i\geqslant 1.\]
The set $\lbrace i\geqslant 1\:\vert\:\alpha_i>1\rbrace$ is infinite and the previous inequality does not occure infinitly, we conclude that there exists $i_0\geqslant 1$ such that $\delta_{i_0}(h)=0$ and so $\mu'_{i_0}(h)=\mu'(h)$.\\\qed
\\ \\\indent From now, we suppose that we have construct an infinite set of key polynomials $\textbf{Q}_{\omega}=\lbrace Q_{i}\rbrace_{i\geqslant 1}$ such that $\alpha_i=1$,for all $i$ sufficiently large. Thus for this $i$, we have:
\[Q_{i+1}=Q_i+z_i,\]
where $z_i$ is an $i$-standard homogeneous expansion, of value $\beta_i$, not containing $Q_i$. More we suppose that $\beta_i\in G_1$ and there exists $h\in K[x]$ such that, for all $i\geqslant 1$: $$\mu'_i(h)<\mu'(h).$$ Write $Q_\omega$ the monic polynomial of smaller degree possible satisfying the previously inequality. As in the proof of the Proposition \ref{sialphainfini}, write:
\[\delta_i(Q_\omega)=\max S_i(Q_\omega,\beta_i),\]
where:
\[S_i(Q_\omega,\beta_i)=\lbrace j\in\lbrace 0,...,s_{i}\rbrace\:\vert\:j\beta_i+\mu'\left(c_{j,i}\right)=\mu'_i(Q_\omega)\rbrace,\]
\[Q_\omega=\sum\limits_{j=0}^{s_{i}}c_{j,i}Q_{i}^j.\]
By (1) of Proposition 5.2 of \cite{spivamahboub}, we have:
\[\alpha_{i+1}\delta_{i+1}(Q_\omega)\leqslant\delta_i(Q_\omega),\:\forall\:i\geqslant 1.\]
Since $\alpha_i=1$ for $i$ sufficiently large, there exists $\delta_\omega\in\mathbb N^*$ such that $\delta_\omega=\delta_i(Q_\omega)$, for $i$ sufficiently large.

\begin{prop}\label{sialphafini}(\cite{spivamahboub}, Proposition 6.8)
Under the previously assumptions, if $\delta_\omega=1$ then the sequence $\lbrace\beta_i\rbrace_{i\geqslant 1}$ is not bounded in $G_1$.
\end{prop}
\section{Defect and effective degree}
In this section, we do not suppose any assumption about the rank of the valuation. We use the construction of key polynomials of W. Mahboub (see \cite{mahboubthese}). Recall that the order type of the set of key polynomials is at most $\omega\times\omega$.
\begin{lem}
Let $(K,\mu)\hookrightarrow (L,\mu')$ be a finite and simple valued field. Let $\{Q_i\}_{i\in\Lambda}$ be the well ordered set of key polynomials associated to this extension. Then there exists $n_0\in\mathbb N^*$ such that $\Lambda\leqslant\omega n_0$, ie $\{Q_i\}_{i\in\Lambda}$ have a finite number of limit key polynomials.
\end{lem}
\noindent\textit{Proof}: Let $P\in K[x]$ monic and irreducible such that $L=K[x]/(P)$. Denote also by $\mu'$ the pseudo-valuation of $K[x]$ for which the set of key polynomials is associated, thus for all $Q\in (P)$, $\mu'(Q)=\infty$. Suppose that the set $\{Q_{\omega i}\}_{i\geqslant 1}$ is infinite. Such $\alpha_{\omega i}>1$, for all $i\geqslant 1$, the sequence $\{d^{\:\circ}(Q_{\omega i})\}_{i\geqslant 1}$ is strictly increasing. Then there exsits $n_0\in\mathbb N^*$ such that $d^{\:\circ}(P)\leqslant d^{\:\circ}(Q_{\omega i})$, for all $i\geqslant n_0$.
\\Let $i\geqslant n_0$, by euclidean division there exists $S_i,R_i\in K[x]$ such that $Q_{\omega i}=PS_i+R_i$ with $R_i=0$ or $0\leqslant d^{\:\circ}(R_i)<d^{\:\circ}(P)\leqslant d^{\:\circ}(Q_{\omega i})$. If $R_i\neq 0$ then $\mu'(Q_{\omega i})=\mu'(R_i)$ which is absurd by minimality of the degree of a limit key polynomial. If $R_i=0$ then $Q_{\omega i}\in (P)$ and the $Q_{\omega i}$ are not key polynomials for $i>n_0$.\\\qed
\begin{lem}\label{lemmegeneraldefaut}
Let $(K,\mu)$ be an henselian field and $L$ be a finite and simple extension of $K$. By definition, $\mu$ extends uniquely to $L$ and this extension corresponds to a (pseudo-)valuation in $K\left[ x\right]$ denoted by $\mu'$. Consider $\lbrace Q_i \rbrace_{i\in\Lambda}$ the well ordered set of key polynomials associated to $\mu'$ and $n_0\in\mathbb N^*$ the smallest possible such that $\Lambda\leqslant \omega n_0$. Then there exists an index $i_0\in\Lambda$ having a predecessor, such that:
\[ \left[ L:K\right]=d^{\:\circ}(Q_{\omega(n_0-1)+i_0})d_{\omega n_0}\]
where:
\[d_{\omega n_0}=\left \{ \begin{array}{ccl}  \delta_{\omega n_0} & \textup{si} & \Lambda=\omega n_0\textup{ et } \sharp\lbrace i\geqslant 1\:\vert\:\alpha_{\omega(n_0-1)+i}=1\rbrace=+\infty \\  1  & \textup{si} & \Lambda<\omega n_0 \textup{ ou } \Lambda=\omega n_0\textup{ et } \sharp\lbrace i\geqslant 1\:\vert\:\alpha_{\omega(n_0-1)+i}=1\rbrace<+\infty \end{array} \right.\] 
\end{lem}
\noindent\textit{Proof}: By assumption, $L=K\left[x\right]/(P(x))$ with $P\in K\left[x\right]$ irreducible and monic.
\\Suppose that $\Lambda$ is an ordinal having an immediate predecessor, to fix ideas denote it by $\omega (n_0-1)+n$, $n\in\mathbb N^*$. By definition, such $\Lambda<\omega n_0$, $d_{\omega n_0}=1$. By construction of key polynomials, $P=Q_{\omega(n_0-1)+n}$. Thus, $i_0=n$ and:
\[\left[ L:K\right]=d^{\:\circ}(P)=d^{\:\circ}(Q_{\omega(n_0-1)+i_0})d_{\omega n_0}.\]
Suppose now that $\Lambda$ is a limit ordinal, denote it by $\omega n_0$. 
\\If $\sharp\lbrace i\geqslant 1\:\vert\:\alpha_{\omega(n_0-1)+i}=1\rbrace<+\infty$, then $d_{\omega n_0}=1$. As in the proof of Proposition \ref{sialphainfini}, there exists $i_0\in\mathbb N^*$ minimal such that $\mu'_{\omega (n_0-1)+i_0}(P)=\mu'(P)$. Thus, if we write the $\omega(n_0-1)+i_0$-standard expansion of $P$, we note that $Q_{\omega(n_0-1)+i_0}$ belongs to the center of the valuation $\mu'$ which is the ideal $(P)$. Note that, by definition and construction of the key polynomials and by choice of $i_0$, $d^{\:\circ}(P)\leqslant d^{\:\circ}(Q_{\omega(n_0-1)+i_0})$. Since $Q_{\omega(n_0-1)+i_0}\in (P)$, we conclude that $Q_{\omega(n_0-1)+i_0}=cP$, $c\in K^*$ and so:
\[\left[ L:K\right]=d^{\:\circ}(P)=d^{\:\circ}(Q_{\omega(n_0-1)+i_0})d_{\omega n_0}.\]
If $\sharp\lbrace i\geqslant 1\:\vert\:\alpha_{\omega(n_0-1)+i}=1\rbrace=+\infty$, then $d_{\omega n_0}=\delta_{\omega n_0}$. In this case, $P=Q_{\omega n_0}$. Take $i_0\in\mathbb N^*$ sufficiently large such that $\alpha_{\omega (n_0-1)+i}=1$ and $\delta_{\omega (n_0-1)+i}(P)=\delta_{\omega (n_0-1)+i+1}(P)$ for all $i\geqslant i_0$. Recall that this common value is denoted by $\delta_{\omega n_0}$. By Proposition 2.12 of \cite{vaquie5}, since $K$ is henselian, the Newton polygon $\Delta_{\omega (n_0-1)+i}(P)$ have a unique side of slope $\beta_{\omega (n_0-1)+i}$. This is equivalent to:
\begin{align*}
\mu'(c_{0,\omega (n_0-1)+i})&=\mu'(c_{s_{\omega (n_0-1)+i},\omega (n_0-1)+i})+s_{\omega (n_0-1)+i}\beta_{\omega (n_0-1)+i}\\&\leqslant \mu'(c_{j,\omega (n_0-1)+i})+j\beta_{\omega (n_0-1)+i},
\end{align*}
with $P=\sum\limits_{j=0}^{s_{\omega (n_0-1)+i}}c_{j,\omega (n_0-1)+i}Q_{\omega (n_0-1)+i}^j$ and $0\leqslant j \leqslant s_{\omega (n_0-1)+i}$. By Corollary 3.23 of \cite{spivamahboub}, $\beta_{\omega (n_0-1)+i}$ determine always a side of $\Delta_{\omega (n_0-1)+i}(P)$ which only have one, so:
\[\mu'_{\omega (n_0-1)+i}(P)=\mu'(c_{0,\omega (n_0-1)+i})=\mu'(c_{s_{\omega (n_0-1)+i},\omega (n_0-1)+i})+s_{\omega (n_0-1)+i}\beta_{\omega (n_0-1)+i}.\]
We deduce, by definition of $\delta_{\omega (n_0-1)+i}(P)$, that:
\[\delta_{\omega (n_0-1)+i}(P)=s_{\omega (n_0-1)+i}.\]
But, for $i\geqslant i_0$, $d^{\:\circ}(Q_{\omega(n_0-1)+i})=d^{\:\circ}(Q_{\omega(n_0-1)+i_0})$ and $\delta_{\omega (n_0-1)+i}(P)=\delta_{\omega n_0}$. Thus:
\[\delta_{\omega n_0}=s_{\omega (n_0-1)+i}=\dfrac{d^{\:\circ}(P)}{d^{\:\circ}(Q_{\omega(n_0-1)+i_0})}.\]\qed
\begin{coro}\label{corogenedefaut}
With the assumptions of Lemma \ref{lemmegeneraldefaut}, we have:
\[\prod\limits_{j=1}^{n_0}d_{\omega j}=d_{L\vert K}(\mu',\mu).\]
\end{coro}
\noindent\textit{Proof}: Recall that we denote by $\alpha_{\omega j}=d^{\:\circ}_{Q_{\omega(j-1)+i_0}}(Q_{\omega j})$ with $i_0$ a $\omega j$-inessentiel index and $j\in\mathbb N^*$. In \cite{vaquie3}, this number corresponds to the jump of order $j$ denote by $s^{(j)}$. By Proposition 3.4.4 of \cite{mahboubthese}, we can always suppose that $\delta_{\omega j}=\alpha_{\omega j}$, for all $j<n_0$. By Proposition 2.9 of \cite{vaquie3}:
\[d^{\:\circ}(Q_{\omega(n_0-1)+i_0})=[\Gamma_{\mu'_{\omega(n_0-1)+i_0}}:\Gamma_\mu][k_{\mu'_{\omega(n_0-1)+i_0}}:k_\mu]\prod\limits_{j=1}^{n_0-1}\alpha_{\omega j}.\]
We verified after that if $\Lambda=\omega(n_0-1)+i_0$ or if $\Lambda=\omega n_0$, we always have:
\[[\Gamma_{\mu'}:\Gamma_\mu]=[\Gamma_{\mu'_{\omega(n_0-1)+i_0}}:\Gamma_\mu],\]
\[[k_{\mu'}:k_\mu]=[k_{\mu'_{\omega(n_0-1)+i_0}}:k_\mu].\]
Thus, applying Lemma \ref{lemmegeneraldefaut}, we obtain:
\[\left[ L:K\right]=\left([\Gamma_{\mu'}:\Gamma_\mu][k_{\mu'}:k_\mu]\prod\limits_{j=1}^{n_0-1}\alpha_{\omega j}\right)d_{\omega n_0}.\]
Since $d_{\omega j}=\delta_{\omega j}=\alpha_{\omega j}$, we concluded that:
\[d_{L\vert K}(\mu',\mu)=\left(\prod\limits_{j=1}^{n_0-1}\alpha_{\omega j}\right)d_{\omega n_0}=\left(\prod\limits_{j=1}^{n_0-1}d_{\omega j}\right)d_{\omega n_0}=\prod\limits_{j=1}^{n_0}d_{\omega j}.\]\qed
\begin{coro}
With the assumptions of Lemma \ref{lemmegeneraldefaut}, we have:
\[\delta_{\omega n_0}=\alpha_{\omega n_0}.\]
\end{coro}
\noindent\textit{Proof}: By Corollary 2.10 of \cite{vaquie3}, the defect is the total jump, applying Corollary \ref{corogenedefaut} we have the equality:
\[\prod\limits_{j=1}^{n_0}\alpha_{\omega j}=d_{L\vert K}(\mu',\mu)=\prod\limits_{j=1}^{n_0}d_{\omega j}.\]
Such $\delta_{\omega j}=\alpha_{\omega j}\neq 0$, for all $j<n_0$, we deduced that:
\[\delta_{\omega n_0}=\alpha_{\omega n_0}.\]\qed
\begin{coro}\label{corodegregene}
Let $(K,\mu)$ be a valued field and $L$ be a finite and simple extension of $K$. Write $\mu^{(1)},...,\mu^{(g)}$ the different extensions of $\mu$ on $L$, its corresponds to a (pseudo-)valuation of $K\left[ x\right]$ denoted by the same way. Consider $\lbrace Q_l^{(i)} \rbrace_{l\in\Lambda^{(i)}}$ the set of key polynomials associated to $\mu^{(i)}$ and $n_0^{(i)}\in\mathbb N^*$ the smallest possible such that $\Lambda^{(i)}\leqslant \omega n_0^{(i)}$, $1\leqslant i\leqslant g$. Then:
\[d_{L\vert K}(\mu^{(i)},\mu)=\prod\limits_{j=1}^{n_0^{(i)}}d_{\omega j}^{(i)}.\]
We deduced that:
\[ \left[ L:K\right]=\sum\limits_{i=1}^ge_if_id_{\omega}^{(i)}d_{\omega 2}^{(i)}...d_{\omega n_0^{(i)}}^{(i)},\]
where $e_i=\left[\Gamma_{\mu^{(i)}}:\Gamma_\mu\right]$ and $f_i=\left[k_{\mu^{(i)}}:k_\mu\right]$.
\end{coro}
\noindent\textit{Proof}: By Proposition \ref{thmdefautdegre}, we know that:
\[\left[L:K\right]=\sum\limits_{i=1}^ge_if_id_i,\]
where $d_i=d_{L\vert K}(\mu^{(i)},\mu)$, $e_i=\left[\Gamma_{\mu^{(i)}}:\Gamma_\mu\right]$ and $f_i=\left[k_{\mu^{(i)}}:k_\mu\right]$. Applying Corollary \ref{corogenedefaut} to the fields $L^{h(\mu^{(i)})}$ and $K^{h(\mu^{(i)})}$ we have the equalities.\\\qed
\begin{coro}
Under the assumptions of Corollary \ref{corodegregene} and denoting by $p=car(k_\mu)$, then, for $1\leqslant i\leqslant g$ and $1\leqslant j\leqslant n_0^{(i)}$, there exists $e_{\omega j}^{(i)}\in\mathbb N$ such that:
$$\delta_{\omega j}^{(i)}=p^{e_{\omega j}^{(i)}}.$$
\end{coro}
\noindent\textit{Proof}: It is an immediate consequence of Corollary \ref{corodegregene} and Remark \ref{ostrow}.\\\qed
\begin{ex}
\textup{We study the example 3.2. of \cite{vaquie3}. Let $k$ be an algebraic closed field of characteristic $0$. Consider $K=k(y)$ equipped with its $y$-adic valuation  denote by $\mu$. The polynomial:
\[P=x^4+y^2x^3+y^3(y^2-2)x^2-y^5x+y^6\]
is irreducible in $K[x]$. The valuation extends only on two valuations $\mu^{(1)}$ and $\mu^{(2)}$ over the field $L=K[x]/(P)$ given by the sequences of key polynomials $\lbrace Q_i^{(1)}\rbrace_{i\in\Lambda^{(1)}}$ and $\lbrace Q_i^{(2)}\rbrace_{i\in\Lambda^{(2)}}$ such that:
\[\begin{array}{lll}Q_1^{(1)}=x &\beta_1^{(1)}=3/2
\\Q_2^{(1)}=x^2-y^3 & \beta_2^{(1)}=7/2
\\Q_i^{(1)}=Q_2^{(1)}+(y^2-u_{i-3})x &\beta_i^{(1)}=3/2+i\:\:;&i\geqslant 3\end{array}\]
and:
\[\begin{array}{lll}Q_1^{(2)}=x &\beta_1^{(2)}=3/2
\\Q_2^{(2)}=x^2-y^3 & \beta_2^{(2)}=9/2
\\Q_i^{(2)}=Q_2^{(2)}+u_{i-2}x &\beta_i^{(2)}=3/2+(i+1)\:\:;&i\geqslant 3\end{array}\]
where $u_0=0$, $u_l=\sum\limits_{j=1}^{l}c_jy^{j+2}$, for all $l\geqslant 1$, $c_1=c_2=1$ and $c_j=4^{j-2}\dfrac{3...(2j-3)}{6...(2j)}$ for all $j\geqslant 3$. Following the notations of the Corollary \ref{corodegregene}, with $\Lambda^{(1)}=\Lambda^{(2)}=\omega$, and $\alpha_i^{(1)}=\alpha_i^{(2)}=1$ for all $i\geqslant 2$, we deduced that: $$d_\omega^{(1)}=d_\omega^{(2)}=1,$$ since $p=1$ if $car(k)=0$.
Moreover, since a field of caracteristic zero is defectless, we found the fact that: $$d_{L\vert K}(\mu^{(1)},\mu)=1=d_\omega^{(1)}$$ and:
$$d_{L\vert K}(\mu^{(2)},\mu)=1=d_\omega^{(2)}.$$
Finally, since $\Gamma_\mu=\mathbb Z$, $\Gamma_{\mu^{(1)}}=\Gamma_{\mu^{(2)}}=(1/2)\mathbb Z$, and $k_\mu=k_{\mu^{(1)}}=k_{\mu^{(2)}}=k$ we deduced that:
\[e_1=e_2=2\textup{ et }f_1=f_2=1.\]
We still find the result of the Corollary \ref{corodegregene}, that is to say:
\[4=[L:K]=e_1f_1d_\omega^{(1)}+e_2f_2d_\omega^{(2)}=2\times1\times1+2\times1\times1.\] }
\end{ex}
\begin{ex}
\textup{We study the example of W. Mahboub given in \cite{mahboub}. Let $k$ be a field of characteristic $p>2$. Denote by $\nu$ the $z$-adic valuation of $k(z)$ and by $\mu$ the valuation of $K=k(z,y)$ define by the set of key polynomials $\{Q_{y,i}\}_{i\geqslant 1}$ given in \cite{mahboub} extending $\nu$. For $e\in\mathbb N^*$, write:
$$f=x^{p^e}-y^2-z.$$
We want to find all the valuations of $L=K[x]/(f)$ whose extends $\mu$. Consider the valuation $\mu'$ of $L$ define by the set of key polynomials $\lbrace Q_i\rbrace_{i\in\Lambda}$:
\[\begin{array}{lll}
Q_1=x&\beta_1=1-\dfrac{1}{4p}
\\Q_i=x^{p^{e-1}}-\sum\limits_{j=0}^ih_j &\beta_i=1-\dfrac{1}{2^{2i}p^i}\:\:;&i\geqslant 2\end{array}\]
where $h_0=0$ and $h_j=\dfrac{Q_{y,2j+1}^2}{z^{2^{2j}p^j-1}}$, for all $j\geqslant 1$. We deduced that, for all $i\geqslant 1$:
\[f=(Q_i+\sum\limits_{j=0}^ih_j)^p-y^2-z=Q_i^p+\sum\limits_{j=0}^ih_j^p-Q_{y,2}=Q_i^p+\dfrac{Q_{y,2i}}{z^{2^{2i-2}p^{i}-p}}.\]
But $\mu(h_i^p)=\mu\left(\dfrac{Q_{y,2i}}{z^{2^{2i-2}p^{i}-p}}\right)$, so:
\[in_i(f)=X^p-in_\mu (h_i^p)=(X-in_\mu(h_i))^p.\]
Thus, for all $i\geqslant 1$, $\delta_i(f)=p$ and so $\delta_\omega=p$. Moreover we are in the situation where $\Lambda=\omega$ and $\alpha_i=1$, for all $i\geqslant 3$. We deduced, with the Corollary \ref{corodegregene}, that:
\[[k_{\mu'}:k_\mu]=p^{e-1};\]
\[d_{L\vert K}(\mu',\mu)=d_\omega=\delta_\omega=p.\]
Since $[L:K]=p^e$, by the Corollary \ref{corodegregene}, we deduced that $\mu'$ is the unique valuation extending $\mu$ on $L$.}
\end{ex}
\begin{ex}
\textup{Let $k$ be a field of characteristic $p>0$. We denote by $K=k\left(\left(y,z\right)\right)$. Consider the valuation $\nu:K^*\rightarrow \mathbb Z^2_{lex}$, trivial on $k$, defined by $\nu(z)=(1,0)$ and $\nu(y)=(0,1)$. Consider the polynomial $P\in K\left[w\right]$ defined by:
\[P(w)=w^p-z^{p^2-p}w+y^pz^{p^2}.\]
Write $L=K\left[w\right]/(P)$, we want to find all the valuations $\nu^{(i)}$ of $L$ enxtending $\nu$ and all the corresponding sets of key polynomial $\lbrace Q_j^{(i)}\rbrace_{j\in\Lambda^{(i)}}$.
\\\indent Find a valuation of $L$ extending $\nu$ is equivalent to find a pseudo-valuation of $K\left[w\right]$ of kernel $P$. Let $\mu$ a such pseudo-valuation. The only possibles values of $w$ are:
\[\mu(w)=\left \{ \begin{array}{ccc}  (p,p)  \\  (p,1) \\ (p,0) \end{array} \right.\]
But if $\mu(w)=(p,1)$ then $\mu(w^p)=\mu(y^pz^{p^2})=(p^2,p)>\mu(z^{p^2-p}w)=(p^2,1)$, the monomial $w^p$ is not of minimal value in $P$. The two only values possible of $w$ are:
\[\mu(w)=\left \{ \begin{array}{ccc}  (p,p)   \\ (p,0) \end{array} \right.\]
~
\\1) \underline{Suppose that $\mu(w)=(p,p)$.}
\\\indent We are in the situation of the example 5.3.2 of \cite{mahboubthese}. In this case, we prove that the key polynomials are:
$$\begin{array}{l}
Q_1^{(1)}=w,\\Q_j^{(1)}=w-z^pu_j,
\end{array}$$
where $u_j=\sum\limits_{i=1}^iy^{p^i}$ and $\mu(Q_j^{(1)})=(p,p^j)$,for all $j\geqslant 2$. Moreover the polynomial $P$ is a limit key polynomial. Finally, the polynomial $in_j(P)$ is irreducible for all $j\geqslant 1$, the valuation extends uniquely since  $\mu(w)=(p,p)$.
\\\indent Thus, a first way to extend $\nu$ over $L$ is to consider the valuation $\nu^{(1)}$ defined by the key polynomials $Q_j^{(1)}$, $j\geqslant 1$ and $Q_\omega^{(1)}=P$.
\\\indent Write $u=\sum\limits_{i\geqslant 1}y^{p^i}$ the limit of the sequence $(u_j)_j$ in $k\left[\left[y,z\right]\right]$. Remark that $u$ is solution of the equation $X-X^p=y^p$. The other solutions of the equation are the $u-l$ with $l\in\lbrace 1,...,p-1\rbrace$. Writing $\alpha_l=z^p(u-l)$, for $l\in\lbrace 0,...,p-1\rbrace$, we can prove that:
\[\prod\limits_{l=0}^{p-1}(w-\alpha_l)=w^p-z^{p^2-p}w+y^pz^{p^2}.\]
So, the roots of $P$ are $\alpha_0,...,\alpha_{p-1}$ and if we write $Q_{\infty}^{(1)}=w-\alpha_0=w-z^pu$ then:
\[P(w)=Q_{\infty}^{(1)}(w)\prod\limits_{l=1}^{p-1}(w-\alpha_l).\]
Finally, since $in_j(P)=-z^{p^2-p}(X-y^{p^i}z^p)$, for all $j\geqslant 1$ and $\Lambda^{(1)}=\omega$, we deduced that: $$d_{L\vert K}(\nu^{(1)},\nu)=d_\omega^{(1)}=\delta_{\omega}^{(1)}=1.$$
~
\\2) \underline{Suppose that $\mu(w)=(p,0)$.}
\\\indent Denote by $Q_1^{(2)}=w$, then:
$$in_1(P)=X\left(X^{p-1}-z^{p^2-p}\right).$$
Write $Q_2^{(2)}=\left(Q_1^{(2)}\right)^{p-1}-z^{p^2-p}$. Then:
$$P=Q_1^{(2)}Q_2^{(2)}+y^pz^{p^2}.$$
The unque choice for the valuation of $Q_2^{(2)}$ is:
\[\mu(Q_2^{(2)})=(p^2-p,p)>(p,0)=\mu(Q_1^{(2)}).\]
In this case, the polynomial $in_2(P)$ is:
\[in_2(P)=in_\mu(Q_1^{(2)})\left(X+\dfrac{in_\mu(y^pz^{p^2})}{in_\mu(Q_1^{(2)})}\right).\]
Remark that $\mu(u)=\mu(y^p)=(0,p)$ then $\mu(\alpha_0)=\mu(z^pu)=(p,p)$. Writing $h=\sum\limits_{m=1}^{p-1}\alpha_0^{p-m}w^{m-1}$, we remark that $h\in K[w]$ represent $\frac{in_\mu(y^pz^{p^2})}{in_\mu(Q_1^{(2)})}$ because $\mu(h)=\mu(w^{p-2}\alpha_0)=(p^2-p,p)$. Then denote by:
\[Q_3^{(2)}=Q_2^{(2)}+h.\]
But, $(w-\alpha_0)h=\alpha_0(w^{p-1}-\alpha_0^{p-1})$. So:
\[(w-\alpha_0)Q_3^{(2)}=w^p-z^{p^2-p}w+\alpha_0z^{p^2-p}-\alpha_0^p.\]
But, $\alpha_0z^{p^2-p}-\alpha_0^p=uz^{p^2}-u^pz^{p^2}=(u-u^p)z^{p^2}=y^pz^{p^2}$ we deduced that:
\[Q_{\infty}^{(1)}Q_3^{(2)}=(w-\alpha_0)Q_3^{(2)}=P.\]
The sequence of key polynomials stop. The second way to extend $\nu$ on $L$ is to consider the valuation $\nu^{(2)}$ defined with the key polynomials $Q_1^{(2)}$, $Q_2^{(2)}$ and $Q_3^{(2)}$. 
\\Finally, since $\Lambda^{(2)}=3<\omega$, we deduced that: $$d_{L\vert K}(\nu^{(2)},\nu)=d_\omega^{(2)}=1.$$
Then we obtain:
\[p=[L:K]=d^{\:\circ}(P)=d^{\:\circ}(Q_{1}^{(1)})d_\omega^{(1)}+d^{\:\circ}(Q_{3}^{(2)})d_\omega^{(2)}=1+(p-1).\]
}
\end{ex}
\begin{ex}
\textup{We will study an example proposed by S.D. Cutkosky in \cite{cutkcontreex}. Let $k$ a field of characteristic $p>0$. Denote by $K=k(u,v)$ and by $L=K[y]/(f)$ where:
\[f=y^{p^2+1}+y^{p^2}-yv^p+u-v^p.\]
We define a valuation $\nu$ over $K$ such that $\nu(u)=p$ with:
$$R_\nu=\bigcup\limits_{i\geqslant 0}A_i,$$
where $A_i=k\left[u_i,v_i\right]_{(u_i,v_i)}$ and:
\[\left \{ \begin{array}{llc}  u_{0}=u \\  v_{0}=v \end{array} \right. ;\:\left \{ \begin{array}{llc}  u_{pi+j}=\dfrac{u_{pi}}{v_{pi+j}^j} \\  v_{pi+j}=v_{pi} \end{array} \right., 1\leqslant j<p\:;\:
\left \{ \begin{array}{llc}  u_{p(i+1)}=v_{pi} \\  v_{p(i+1)}=\dfrac{u_{pi}}{v_{pi}^p}-\gamma_{p(i+1)} \end{array}\right. ;\: \gamma_{p(i+1)}\in k^*.\]
Then we have $\nu(v_{pi})=\dfrac{1}{p^i}$ and $\nu(u_{pi})=\dfrac{1}{p^{i-1}}$. Remark that, for $i\geqslant 0$:
$$v_{p(i-1)}-\gamma_{p(i+1)}v_{pi}^p=v_{pi}^pv_{p(i+1)},$$
where we write $v_{-p}=u$. Then:
\[\nu(v_{p(i-1)}-\gamma_{p(i+1)}v_{pi}^p)=\dfrac{1}{p^{i-1}}-\dfrac{1}{p^{i+1}}.\]
We want to extend $\nu$ on $L$. There exist a first valuation denoted by $\nu^*$ defined by the sequence of key polynomials $\{Q_l\}_{l\in\Lambda}$:
$$Q_1=y,\:\beta_1=\nu^*(Q_1)=\dfrac{1}{p};$$
$$Q_l=Q_{l-1}-\gamma_{p+2l-3}^{\frac{1}{p}}\prod\limits_{j=0}^{l-2}v_{p(2j+1)},\:\beta_l=\nu^*(Q_l)=\sum\limits_{j=0}^{l-1}\dfrac{1}{p^{2j+1}}.$$
Remark that $\lim\limits_{l\rightarrow +\infty}\beta_l=\dfrac{p}{p^2-1}$, so there exists a first limit key polynomial which is:
$$Q_\omega=y^p-v,\: \beta_\omega=\nu^*(Q_\omega)=1.$$
We continue to define the set of key polynomials, for $i\geqslant 1$, by:
$$Q_{\omega+i}=Q_{\omega+i-1}+\gamma_{p+2i-2}^{\frac{1}{p}}\prod\limits_{j=0}^{i-1}v_{p(2j)},\:\beta_{\omega+i}=\nu^*(Q_{\omega+i})=\sum\limits_{j=0}^i\dfrac{1}{p^{2j}}.$$
Remark that $\lim\limits_{i\rightarrow +\infty}\beta_{\omega +i}=1+\dfrac{1}{p^2-1}$, so there exists a second limit key polynomial which is:
$$Q_{\omega 2}=y^{p^2}+u-v^p=Q_\omega^p+u,\: \beta_{\omega 2}=\nu^*(Q_{\omega 2})=p+\dfrac{1}{p}.$$
We end the set of key polynomials with, for $n\geqslant 1$:
$$Q_{\omega 2+n}=Q_{\omega 2+n-1}+(-1)^nuy^n,\:\beta_{\omega 2+n}=\nu^*(Q_{\omega 2+n})=p+\dfrac{n+1}{p}.$$
The last and third limit key polynomial is $Q_{\omega 3}=f$ and $\beta_{\omega 3}=+\infty$. Remark that for all $n\geqslant 0$:
$$f=(y+1)Q_{\omega 2+n}+(-1)^{n+1}uy^{n+1}.$$
So, in $K\left[\left[y\right]\right]$, $f=(y+1)Q_\infty$ where $Q_\infty=y^{p^2}+u-v^p-uy\sum\limits_{n\geqslant 1}(-1)^ny^n$. Remark that for $l\geqslant 1$:
\[y^p-v=Q_l^p-v_{p(2(l-1))}\prod\limits_{j=0}^{l-2}v_{p(2j+1)}^p,\]
then:
\[in_l(y^p-v)=X^p-in_{\nu}\left(v_{p(2(l-1))}\prod\limits_{j=0}^{l-2}v_{p(2j+1)}^p\right)=\left(X-in_\nu\left(\prod\limits_{j=0}^{l-2}v_{p(2j+1)}\right)\right)^p.\]
We deduce that for all $l\geqslant 1$:
\[\alpha_l=1,\]
\[d_{\omega}=\delta_\omega=\delta_l(y^p-v)=p=\alpha_\omega.\]
In the same way, for $i\geqslant 1$:
\[y^{p^2}+u-v^p=Q_\omega^p+u=Q_{\omega+i}^p-v_{p(2i-1)}\prod\limits_{j=0}^{i-2}v_{p(2j)}^p,\]
\[in_{\omega+i}(y^{p^2}+u-v^p)=X^p-in_{\nu}\left(v_{p(2i-1)}\prod\limits_{j=0}^{i-2}v_{p(2j)}^p\right)=\left(X-in_\nu\left(\prod\limits_{j=0}^{i-1}v_{p(2j)}\right)\right)^p.\]
We deduce that for all $i\geqslant 1$:
\[\alpha_{\omega+i}=1,\]
\[d_{\omega 2}=\delta_{\omega 2}=\delta_{\omega 2+i}(y^{p^2}+u-v^p)=p=\alpha_{\omega 2}.\]
Finally, for all $n\geqslant 0$:
\[in_{\omega 2+n}(f)=X+in_\nu\left((-1)^{n+1}uy^{n+1}\right).\]
We deduce that for all $n\geqslant 1$:
\[\alpha_{\omega 2+n}=1,\]
\[d_{\omega 3}=\delta_{\omega 3}=\delta_{\omega 3+n}(f)=1.\]
By the Corollary \ref{corogenedefaut}:
\[d_{L\vert K}(\nu^*,\nu)=d_{\omega} d_{\omega 2}d_{\omega 3} =p\times p\times 1 =p^2.\]
Remark that $[L:K]=p^2+1$, so there exists a second valuation $\nu^{(2)}$ of $L$ who extends $\nu$ and defectless. By the equality $f=(y+1)Q_\infty$ in $K\left[\left[y\right]\right]$, the valuation verifies $\nu^{(2)}(y)=0$, so this is the trivial valuation. }
\end{ex}
\section{Key polynomials and defectless fields}
Consider $K\hookrightarrow K(x)$ a simple transcendental field extension. Let $\mu'$ be a valuation of $K(x)$, write  $\mu:=\mu'_{\vert\:K}$. We denote by $G$ the value group of $\mu'$ and $G_{1}$ the smallest isolated non zero subgroup of $G$. We suppose that $\mu$ is of rank $1$, $\mu'(x)>0$.
\begin{prop}\label{polycleencar0}
If $K$ is defectless, there exists a $1$-complete set of key polynomials $\lbrace Q_i\rbrace_{i\in\Lambda}$ such that $\Lambda$ is a finite set or $\mathbb N^*$. In particular, there is no limit key polynomial  for valuations of rank $1$ over defectless fields.
\end{prop}
\noindent\textit{Proof}: Apply the process of \cite{spivamahboub}. If there exists $i_0\in\mathbb N$, such that $\beta_{i_0}\notin G_1$, we write $\Lambda=\lbrace 1,...,i_0-1\rbrace$ and, by definition, $\lbrace Q_i\rbrace_{i\in\Lambda}$ is $1$-completed. Otherwise, for all $i\in\mathbb N$, $\beta_i\in G_1$ and we write $\Lambda=\mathbb N^*$. If $\sharp\lbrace i\geqslant 1\:\vert\:\alpha_i>1\rbrace=+\infty$, by Proposition \ref{sialphainfini}, the set $\lbrace Q_i\rbrace_{i\in\Lambda}$ is $1$-completed. If $\sharp\lbrace i\geqslant 1\:\vert\:\alpha_i>1\rbrace<+\infty$, write $Q_\omega$ the monic polynomial of smallest degree such that, for all $i\geqslant 1$: $$\mu'_i(Q_\omega)<\mu'(Q_\omega).$$
Since $K$ is defectless, the extension $K\hookrightarrow L=K[x]/(Q_\omega)$ is defectless. Consider $\mu'$ as the composition of a valuation $\mu^{(1)}$ of value group $G_1$ centered on $K[x]/(Q_\omega)$ and a valuation $\theta$ of value group $G/G_1$ centered on $K[x]_{(Q_\omega)}$. The set of key polynomials for $\mu^{(1)}$ is the same as the set of key polynomials for $\mu'$ except that $\mu^{(1)}(Q_\omega)=\infty$. Thus, $d_{L\vert K}(\mu^{(1)},\mu)=1$. By Corollary \ref{corodegregene}, we deduce that: $$\delta_\omega=d_{L\vert K}(\mu^{(1)},\mu)=1.$$
Applying Proposition \ref{sialphafini}, we conclude that the sequence $\lbrace\beta_i\rbrace_{i\geqslant 1}$ is unbounded in $G_1$ and so, by Proposition \ref{sibornealorscomplet}, the set $\lbrace Q_i\rbrace_{i\in\Lambda}$ is a $1$-completed set of key polynomials.\\\qed
\begin{coro}
If $car(k_\mu)=0$, there exists a $1$-completed set of key polynomials $\lbrace Q_i\rbrace_{i\in\Lambda}$ such that $\Lambda$ is a finite or $\mathbb N^*$. In particular, there is no limit key polynomials for valuations of rank $1$ whose residual field is of characteristic zero.
\end{coro}
\noindent\textit{Proof}: Apply Proposition \ref{cropasdefaut} and Proposition \ref{polycleencar0}.\\\qed

\section{Local uniformization of quasi-excellent local domain without defect}
\indent We extend here the results of the section 8 of \cite{jccar0} for a valuation satisfying some inductive assumption about defect. More precisely, in \cite{jccar0}, to obtain a theorem of monomialization, the valuation needs to have a complete set of key polynomials without limit key polynomial: this is the case if the valuation is defectless. As a  corollary, we find the local uniformization in characteristic zero.
\\ \\\indent Let $(R,\mathfrak{m},k)$ be a local complete regular equicharacteristic ring of dimension $n$ with $\mathfrak{m}=\left(u_1,...,u_{n}\right)$. Let $\nu$ be a valuation of $K=Frac(R)$, centered on $R$, of value group $\Gamma$ and $\Gamma_1$ the smallest non-zero isolated subgroup of $\Gamma$. Write: 
\[H=\lbrace f\in R\:\vert\: \nu(f)\notin\Gamma_{1}\rbrace.\]
$H$ is a prime ideal of $R$ (see the proof of Theorem \ref{thmprelimcar0}). Moreover suppose that:
\[n=e(R,\nu)=emb.dim\left(R/H\right),\]
 that is to say:
 \[H\subset\mathfrak{m}^2.\]
Write $r=r(R,u,\nu)=\dim_\mathbb Q\left(\sum\limits_{i=1}^n\mathbb Q\nu(u_i)\right)$.
\\ The valuation $\nu$ is unique if $ht(H)=1$; it is the composition of the valuation $\mu:L^{*}\rightarrow\Gamma_{1}$ of rank $1$ centered on $R/H$, where $L=Frac(R/H)$, with the valuation $\theta :K^{*}\rightarrow \Gamma / \Gamma_{1}$, centered on $R_{H}$, such that $k_{\theta}\simeq \kappa(H)$.
\\By abuse of notation, for $f\in R$, we will denote by $\mu(f)$ instead of $\mu(f\mod H)$.
By the Cohen's theorem, we can suppose that $R$ is of the form:
\[R=k\left[ \left[ u_{1},...,u_{n}\right] \right].\]
For $j\in \lbrace r+1,...,n\rbrace$, write $\lbrace Q_{j,i}\rbrace_{i\in\Lambda_{j}}$ the set of key polynomials of the extension $k\left( \left( u_{1},...,u_{j-1}\right) \right)\hookrightarrow k\left( \left( u_{1},...,u_{j-1}\right) \right)(u_{j})$, $\textbf{Q}_{j,i}=\left\lbrace Q_{j,i'}\vert i'\in\Lambda_{j},i'<i\right\rbrace $, $\Gamma^{(j)}$ the value group of $\nu_{\vert k\left( \left( u_{1},...,u_{j}\right) \right)}$ and $\nu_{j,i}$ the $i$-troncation of $\nu$ for this extension.
\\ \\\indent For the definition of local framed sequences, one may consult D\'efinition 7.1 and the sections 4.1 and 4.2 of \cite{jccar0}.
\begin{thm}\label{thmeclatformcar0}
Suppose that, for $R_{n-1}=k\left[ \left[ u_{1},...,u_{n-1}\right] \right]$ we have:
\begin{enumerate}
\item 
\begin{enumerate}
\item Or $H\cap R_{n-1}\neq (0)$ and there exists a local framed sequence $(R_{n-1},u)\rightarrow (R',u')$ such that: $$e(R',\nu)<e(R_{n-1},\nu);$$
\item Or $H\cap R_{n-1}=(0)$ and for all $f\in R_{n-1}$, there exists a local framed sequence $(R_{n-1},u)\rightarrow (R',u')$ such that $f$ is a monomial in $u'$ times a unit of $R'$.
\end{enumerate}
\item The local framed sequence $(R_{n-1},u)\rightarrow (R',u')$ of (1) can be chosen defined over $T$.
\end{enumerate}
Moreover suppose that $k\left( \left( u_{1},...,u_{n-1}\right) \right)\hookrightarrow k\left( \left( u_{1},...,u_{n}\right) \right)/H$ is defectless. Then the assumptions 1. and 2. are true with $R$ instead of $R_{n-1}$.
\end{thm}
\noindent\textit{Proof}: The proof is the same as the proofs on Theorem 5.1 and 7.2 of \cite{jccar0}. With the assumptions of the Theorem \ref{thmeclatformcar0}, we can use the Proposition 5.2 of \cite{jccar0}. Then $H$ is generated by a irreducible monic polynomial in $u_n$. Since $k\left( \left( u_{1},...,u_{n-1}\right) \right)\hookrightarrow k\left( \left( u_{1},...,u_{n}\right) \right)/H$ is defectless, by Proposition \ref{polycleencar0}, the set of key polynomials $\lbrace Q_{j,i}\rbrace_{i\in\Lambda_{j}}$ has not limit key polynomial. To conclude it is sufficient to apply Theorem 7.2 of \cite{jccar0}.\\\qed
\\\indent Consequently, we obtain the local uniformization of a valuation of rank $1$ centered on a local quasi-excellent equicharacteristic domain satisfaying some assumptions of lack of defect. The proof uses the notion of implicit prime ideal, for more details see \cite{ideimp} or section 4.3 of \cite{jccar0}.
\begin{defi}
For a local noetherian ring $(R,\mathfrak{m})$, with $\mathfrak m=(u)=(u_1,...,u_n)$ and $f_1,...,f_s\in\mathfrak m$, we call the \textbf{monomial property for R and }$\boldsymbol{f_1,...,f_s}$ this three statements:
\begin{enumerate}
\item $R$ is regular;
\item For $1\leqslant j\leqslant s$, $f_j$ is a monomial in $u$ times a unit of $R$;
\item For $1\leqslant j\leqslant s$, $f_1$ divides $f_j$ in $R$.
\end{enumerate}
\end{defi}
\begin{defi}
Let $(S,\mathfrak m,k)$ be a local domain of quotient field $L$ and $\mu$ be a valuation of $L$ of rank $1$ and of value group $\Gamma_1$, centered on $S$.
\\Write $u=(u_1,...,u_n)$ a  minimal set of generators of $\mathfrak m$ and $\overline H$ the implicit prime ideal of $\widehat S$.
\\Write $u=(y,x)$ such that $x=(x_1,...,x_l)$, $l=emb.dim\left(\widehat{S}/\overline{H}\right)$ and the images of $x_1,...,x_l$ in $\widehat S/\overline H$ induce a minimal set of generators of $(\mathfrak m \widehat S)/\overline H$. 
\\Let $f_1,...,f_s\in\mathfrak m$ such that $\mu(f_1)=\min_{1\leqslant j\leqslant s}\lbrace\mu(f_j)\rbrace$. 
\begin{enumerate}
\item We say that $\boldsymbol S$ \textbf{and} $\boldsymbol{f_1,...,f_s}$ \textbf{have the quotient local uniformization property of dimension l} if there exists a local framed sequence:
\[ \xymatrix{\left( S,u,k\right)=\left( S^{(0)},u^{(0)},k^{(0)}\right) \ar[r]^-{\rho_{0}} & \left( S^{(1)},u^{(1)},k^{(1)}\right) \ar[r]^-{\rho_{1}} & \ldots   \ar[r]^-{\rho_{i-1}} & \left( S^{(i)},u^{(i)},k^{(i)}\right)}, \]
such that $\widehat{S^{(i)}}/\overline H^{(i)}$ and $\overline f_1,...,\overline f_s$ have the monomial property, where $\overline H^{(i)}$ is the implicit prime ideal of $\widehat{S^{(i)}}$ and $\overline{f_j}$ are the images of $f_j\mod\overline H^{(i)}$, $1\leqslant j\leqslant s$.
\item We say that $\boldsymbol S$ \textbf{and} $\boldsymbol{f_1,...,f_s}$ \textbf{have the local uniformization property of dimension l} if there exists a local framed sequence:
\[ \xymatrix{\left( S,u,k\right)=\left( S^{(0)},u^{(0)},k^{(0)}\right) \ar[r]^-{\rho_{0}} & \left( S^{(1)},u^{(1)},k^{(1)}\right) \ar[r]^-{\rho_{1}} & \ldots   \ar[r]^-{\rho_{i-1}} & \left( S^{(i)},u^{(i)},k^{(i)}\right)}, \]
such that $S^{(i)}$ and $\overline f_1,...,\overline f_s$ have the monomial property.
\end{enumerate}
\end{defi}
\begin{thm}\label{thmprelimcar0}
Let $(S,\mathfrak m,k)$ be a local noetherian equicharacteristic domain of quotient field $L$ and $\mu$ be a valuation of $L$ of rank $1$ with value group $\Gamma_1$, centered on $S$.
\\Write $u=(u_1,...,u_n)$ a minimal set of generators of $\mathfrak m$ and $\overline H$  the implicit prime ideal of $\widehat S$.
\\Write $u=(y,x)$ such that $x=(x_1,...,x_l)$, $l=emb.dim\left(\widehat{S}/\overline{H}\right)$ and the images of $x_1,...,x_l$ in $\widehat S/\overline H$ induce a minimal set of generators of $(\mathfrak m \widehat S)/\overline H$. 
\\Suppose that the quotient local uniformization property of dimension l-1 is true for all local domain and $k\left(\left(x_1,...,x_{l-1}\right)\right)\hookrightarrow k\left(\left(x_1,...,x_{l-1}\right)\right)[x_l]/\overline H$ is defectless.
\\Let $f_1,...,f_s\in\mathfrak m$ such that $\mu(f_1)=\min_{1\leqslant j\leqslant s}\lbrace\mu(f_j)\rbrace$. Then $S$ and $f_1,...,f_s$ have the quotient local uniformization property of dimension l.

\end{thm}
\noindent\textit{Proof}: By Theorem 2.1 of \cite{ideimp}, $\mu$ extends uniquely in a valuation $\widehat\mu$ centered on $\widehat S/\overline H$. 
\\By the Cohen's structure theorem, we know that there exists a complete regular local ring of characteristic zero $R$ and a surjective morphism $\varphi$:
\[\varphi:R\twoheadrightarrow \widehat S/\overline H.\]
Write $H=\ker \varphi$, since $\overline H$ is a prime ideal (Theorem 2.1 of \cite{ideimp}), $H$ is a prime ideal of $R$. Choose $R$ such that $\dim (R)=l$. Write $K$ the quotient field of $R$, $K$ is of the form $k\left(\left(x_1,...,x_{l}\right)\right)$.Let $\theta$ a valuation of $K$, centered on $R_H$, such that $k_\theta=\kappa(H)$. If we look $\widehat\mu$ as a valuation centered on $R/H$ by the morphism $\varphi$, we can consider the valuation $\nu=\widehat{\mu}\circ\theta$ centered on $R$ and of value group $\Gamma$. Then, $\Gamma_1$ is the non-zero smallest isolated subgroup of $\Gamma$ and:
\[H=\lbrace f\in R\:\vert\:\nu(f)\not\in\Gamma_1\rbrace.\]
\\Let $T=\varphi^{-1}(\sigma(S))$, this is a local subring $R$ of maximal ideal $\varphi^{-1}(\sigma(\mathfrak m))=\mathfrak m\cap T$. Then, $T$ contains $x_1,...,x_l$ and:
\[T/(\mathfrak m\cap T)\simeq k.\]
Since $k\left(\left(x_1,...,x_{l-1}\right)\right)\hookrightarrow k\left(\left(x_1,...,x_{l-1}\right)\right)[x_l]/ H$, we can apply Theorem \ref{thmeclatformcar0}. We end the proof in the same way that the proof of Theorem 8.1 of \cite{jccar0}.\\\qed
\begin{thm}\label{uniflocalerang1car0}
Let $(S,\mathfrak m,k)$ be an equicharacteristic quasi-excellent local domain of quotient field $L$ and $\mu$ be a valuation of $L$ of rank $1$ and of value group $\Gamma_1$, centered on $S$.
\\Write $u=(u_1,...,u_n)$ a  minimal set of generators of $\mathfrak m$ and $\overline H$ the implicit prime ideal of $\widehat S$.
\\Write $u=(y,x)$ such that $x=(x_1,...,x_l)$, $l=emb.dim\left(\widehat{S}/\overline{H}\right)$ and the images of $x_1,...,x_l$ in $\widehat S/\overline H$ induce a minimal set of generators of $(\mathfrak m \widehat S)/\overline H$. 
\\Suppose that the local uniformization property of dimension l-1 is true for all local domain and $k\left(\left(x_1,...,x_{l-1}\right)\right)\hookrightarrow k\left(\left(x_1,...,x_{l-1}\right)\right)[x_l]/\overline H$ is defectless.
\\Let $f_1,...,f_s\in\mathfrak m$ such that $\mu(f_1)=\min_{1\leqslant j\leqslant s}\lbrace\mu(f_j)\rbrace$. Then $S$ and $f_1,...,f_s$ have the local uniformization property of dimension l.
\\In other words, $\mu$ admits an embedded local uniformization in the sense of Property 2.11 of \cite{jccar0}.
\end{thm}
\noindent\textit{Proof}: Take the notations of Theorem \ref{thmprelimcar0}. We saw that there exists a surjective morphism:
\[\psi:\widehat S\twoheadrightarrow \widehat{S}/\overline{H}\simeq R/H.\]
Since $k\left(\left(x_1,...,x_{l-1}\right)\right)\hookrightarrow k\left(\left(x_1,...,x_{l-1}\right)\right)[x_l]/\overline H$ is defectless, by Theorem \ref{thmprelimcar0}, after an auxiliary local framed sequence, we can suppose that $\widehat{S}/\overline{H}$ is regular and so that $R/H\simeq k\left[\left[x_{1},...,x_l\right]\right]$. The end of the proof is the same 
as the proof of Theorem 8.3 of \cite{jccar0}.\\\qed
\begin{coro}
Let $(S,\mathfrak m,k)$ be a quasi-excellent local domain of quotient field $L$ and $\mu$ be a valuation of $L$ of rank $1$ centered on $S$, such that $car\left(k_\mu\right)=0$.
\\Then $\mu$ admits an embedded local uniformization in the sense of Property 2.11 of \cite{jccar0}.
\end{coro}
\noindent\textit{Proof}: Write $u=(u_1,...,u_n)$ a minimal set of generators of $\mathfrak m$ and $\overline H$ the implicit prime ideal of $\widehat S$. Write $u=(y,x)$ such that $x=(x_1,...,x_l)$, $l=emb.dim\left(\widehat{S}/\overline{H}\right)$ and the images of $x_1,...,x_l$ in $\widehat S/\overline H$ induce a minimal set of generators of $(\mathfrak m \widehat S)/\overline H$. By Theorem 2.1 of \cite{ideimp}, $\mu$ etends uniquely in a valuation $\widehat\mu$ centered on $\widehat S/\overline H$ and, since $car\left(k_\mu\right)=0$, then $car(k_{\:\widehat\mu})=0$. We saw that there exists a surjective morphism:
\[\psi:\widehat S\twoheadrightarrow \widehat{S}/\overline{H}\simeq R/H,\]
where $H=\ker \psi$. The quotient field of $R$ is of the form $k\left(\left(x_1,...,x_{l}\right)\right)$. By Proposition \ref{cropasdefaut}, we deduce that the valued fields $\left(k\left(\left(x_1,...,x_{j-1}\right)\right),\widehat\mu_{\vert k\left(\left(x_1,...,x_{j-1}\right)\right)}\right)$ are defectless. To conclude, it is sufficents to apply recursively on $j\in \{1,...,n \}$ the Theorem \ref{uniflocalerang1car0} for $\left(k\left(\left(x_1,...,x_{j-1}\right)\right),\widehat\mu_{\vert k\left(\left(x_1,...,x_{j-1}\right)\right)}\right)$.\\\qed
\bibliographystyle{plain}
\bibliography{biblio}

\begin{thebibliography}{10}

\bibitem{cutkcontreex}
Steven~Dale Cutkosky.
\newblock \textit{Ramification of valuations and counterexamples to local
  monomialization in positive characteristic}, 2014.
\newblock preprint math.AG/arXiv:1404.7459.

\bibitem{favrejon}
Charles Favre and Mattias Jonsson.
\newblock {\em \textit{The valuative tree}}, volume 1853 of {\em
  \textup{Lecture Notes in Mathematics}}.
\newblock Springer-Verlag, Berlin, 2004.

\bibitem{spivamahboub}
F.~J. Herrera~Govantes, W.~Mahboub, M.~A. Olalla~Acosta, and M.~Spivakovsky.
\newblock \textit{Key polynomials for simple extensions of valued fields},
  2014.
\newblock preprint math.AG/arXiv:1406.0657.

\bibitem{spivaherrera}
F.~J. Herrera~Govantes, M.~A. Olalla~Acosta, and M.~Spivakovsky.
\newblock \textit{Valuations in algebraic field extensions}.
\newblock {\em \textup{J. Algebra}}, 312(2):1033--1074, 2007.

\bibitem{ideimp}
F.~J. Herrera~Govantes, M.~A. Olalla~Acosta, M.~Spivakovsky, and B.~Teissier.
\newblock \textit{Extending a valuation centered in a local domain to the
  formal completion}.
\newblock {\em \textup{Proc. London Math. Soc.}}, 105(3):571--621, 2012.

\bibitem{kulhlivre}
Franz-Viktor Kuhlmann.
\newblock {\em \textit{Valuation Theory}}.
\newblock \textup{http://math.usask.ca/~fvk/Fvkbook.htm}, 2011.

\bibitem{mahboub}
Wael Mahboub.
\newblock \textit{Key Polynomials}.
\newblock {\em \textup{J. Pure Appl. Algebra}}, 217(6):989--1006, 2013.

\bibitem{mahboubthese}
Wael Mahboub.
\newblock {\em \textit{Une construction explicite de polyn\^omes-cl\'es pour
  des valuations de rang fini}}.
\newblock \textup{Th\`ese de Doctorat, Institut de Math\'ematiques de
  Toulouse}, 2013.

\bibitem{jccar0}
Jean-Christophe San~Saturnino.
\newblock \textit{Uniformisation locale des sch\'emas quasi-excellents de
  caract\'eristique nulle}.
\newblock {\em \textup{arXiv:1311.3525}}, 2013. Soumis.

\bibitem{vaquie2}
Michel Vaqui{\'e}.
\newblock \textit{Extension d'une valuation}.
\newblock {\em \textup{Trans. Amer. Math. Soc.}}, 359(7):3439--3481
  (electronic), 2007.

\bibitem{vaquie3}
Michel Vaqui{\'e}.
\newblock \textit{Famille admissible de valuations et d\'efaut d'une
  extension}.
\newblock {\em \textup{J. Algebra}}, 311(2):859--876, 2007.

\bibitem{vaquie5}
Michel Vaqui{\'e}.
\newblock \textit{Extensions de valuation et polygone de {N}ewton}.
\newblock {\em \textup{Ann. Inst. Fourier (Grenoble)}}, 58(7):2503--2541, 2008.

\end{thebibliography}
\end{document}